\documentclass{elsarticle-modified}
\usepackage[usenames,dvipsnames]{color}
\usepackage{graphicx}
\usepackage{amsmath}
\usepackage{eufrak}
\usepackage{ulem}
\usepackage{tikz}
\usepackage{epstopdf}
\newcommand\rev[1]{{#1}}

\newtheorem{remark}{Remark}

\newcommand\Tstrut{\rule{0pt}{2.6ex}}         
\newcommand\Bstrut{\rule[-0.9ex]{0pt}{0pt}}   

\usepackage{tikz,pgfplots}
        \pgfplotsset{compat = 1.3}
        \pgfplotsset{minor grid style={dotted}} \pgfplotsset{major grid
        style={dashed}}

        \pgfplotsset{every x tick label/.append style={font=\footnotesize,
        yshift=0.25ex}}

        \pgfplotsset{every y tick label/.append
        style={font=\footnotesize, xshift=0.25ex}}

        \definecolor{colorclassyorange}{rgb}{0.95000,0.32500,0.09800}
        \definecolor{colorchromeyellow}{rgb}{1.00000,0.6549,0}%
        \definecolor{colorpaleyellow}{rgb}{1.00000,0.8549,0.1}%

        \definecolor{colorclassyblue}{rgb}{0.00000,0.44706,0.74118}%
        \definecolor{colorpurple}{rgb}{0.49400,0.18400,0.55600}%
        \definecolor{colorfuschia}{rgb}{0.95039,0.0,0.95039}%
        \definecolor{colorlemongreen}{rgb}{0.6,0.8,0}%
        \definecolor{colorreal}{rgb}{0.92941,0.79412,0.12549}%
        \definecolor{colorimag}{rgb}{0.00000,0.49804,0.00000}%
        \definecolor{colorabs}{rgb}{1.00000,0.00000,0.00000}%

\newcommand\dd{\mathrm{d}} 
\font\Bbb=msbm10 scaled 1000
\def\BB#1{\mbox{\Bbb #1}}              

\newcommand{\ee}{{\mathrm e}}
\newcommand{\ii}{{\mathrm i}}

\newcommand{\K}{\mathcal{K}}

\newcommand{\DDD}{\mathcal{D}}
\newcommand{\WWW}{\mathcal{W}}
\newcommand{\SSS}{\mathcal{S}}
\newcommand{\EEE}{\mathcal{E}}
\newcommand{\LLL}{\mathcal{L}}
\newcommand{\dx}{\partial_x}

\def\O#1{{\mathcal O}\!\left(#1\right)}

\newcommand{\ang}[2]{\left\langle #1 \right\rangle_{#2}}

\newcommand{\ve}{\varepsilon}
\newcommand{\schr}{Schr\"odinger }

\newcommand{\RRR}{\mathcal{R}}

\newcommand{\Hc}[1]{#1}

\newcommand{\dt}{\Hc{h}}

\newcommand\nomark[1]{#1}

\newcommand{\norm}[2]{\left\| #2 \right\|_{#1}}

\numberwithin{equation}{section}

\def\schr{Schr{\"o}dinger\ }

\journal{Appl. Math. Comput.}

\begin{document}

\begin{frontmatter}

\title{Time adaptive Zassenhaus splittings
for the Schr{\"o}dinger equation in the semiclassical regime\tnoteref{t1}}

\tnotetext[t1]{This work was supported in part by
the Vienna Science and Technology Fund (WWTF) under grant MA14-002.
The work of Karolina Kropielnicka on this project
was financed by The National Center for Science,
under grant no. 2016/23/D/ST1/02061.}

\author{Winfried Auzinger\corref{ourcorrespondingauthor}}
\address{Technische Universit{\"a}t Wien, Institut f{\"u}r Analysis und Scientific Computing,
Wiedner Hauptstrasse 8--10/E101, A-1040 Wien, Austria}
\ead{w.auzinger@tuwien.ac.at}
\ead[url]{http://www.asc.tuwien.ac.at/ winfried}

\author{Harald Hofst{\"a}tter\corref{}}
\address{Universit{\"a}t Wien, Institut f{\"u}r Mathematik, Oskar-Morgenstern-Platz 1, A-1090 Wien, Austria}
\ead{hofi@harald-hofstaetter.at}
\ead[url]{http://www.harald-hofstaetter.at}

\author{Othmar Koch\corref{}}
\address{Universit{\"a}t Wien, Institut f{\"u}r Mathematik, Oskar-Morgenstern-Platz 1, A-1090 Wien, Austria}
\ead{othmar@othmar-koch.org}

\author{Karolina Kropielnicka}
\address{Institute of Mathematics, University of Gdansk, ul. Wit Stwosz 57, 80-308 Gdansk}
\ead[url]{https://mat.ug.edu.pl/ kmalina}
\ead{kmalina@mat.ug.edu.pl}

\author{Pranav Singh}
\address{Mathematical Institute, University of Oxford, Andrew Wiles Building,
Radcliffe Observatory Quarter, Woodstock Road, Oxford OX2 6GG}
\ead{pranav.singh@trinity.ox.ac.uk}
\ead[url]{https://www.maths.ox.ac.uk/people/pranav.singh}

\cortext[ourcorrespondingauthor]{Corresponding author}

\begin{abstract}
Time dependent \schr equations with \rev{conservative force field} commonly
constitute a major challenge in the numerical approximation, especially
when they are analysed in the semiclassical regime. Extremely high
oscillations originate from the semiclassical parameter, and call for
appropriate methods. We propose to employ a combination of asymptotic
Zassenhaus splitting with time adaptivity. While the former turns the
disadvantage of the semiclassical parameter into an advantage, leading to
highly efficient methods with low error constants, the latter enables to
choose an optimal time step and to speed up the calculations when the
oscillations subside. We support the results with numerical examples.
\end{abstract}



\end{frontmatter}

\section{Introduction} \label{sec:intro}

In this paper we are concerned with developing a time-adaptive method for
solving the \schr equation in the semiclassical regime,
\begin{equation} \label{eq:TDSE}
\begin{aligned}
\partial_t \psi(x,t) &=
\ii\,\ve \dx^2 \psi(x,t) - \ii\,\ve^{-1} V(x) \psi(x,t)\, \rev{=: H\psi(x,t)},\quad x\in I \subseteq \BB{R},\ t\geq 0, \\
\psi(x,0) &= \psi_0(x),
\end{aligned}
\end{equation}
with periodic boundary conditions imposed on the interval $I \subseteq
\BB{R}$. The interaction potential $V(x)$ is a real and periodic function.
The regularity required of $V$ and $\psi$ depends on the desired order of
the numerical method. For the sake of simplicity, we assume $V\in \mathrm{C}_{\mathrm{p}}^{\infty}(I;\BB{R})$ and $\psi\in
\mathrm{C}_{\mathrm{p}}^{\infty}(I;\BB{C})$, where the  subscript ``$\mathrm{p}$'' denotes periodicity.

The semiclassical parameter $\ve$ induces oscillations of wavelength $\O{\ve}$ in
the solution $\psi$, both in space and in time
\cite{bao02ots,jin11mcm,Singh17}. In this regime, finite difference schemes
require an excessively fine spatial grid and very small time steps
\cite{Markowich1999}. Consequently, they are found to be ineffective in
comparison with spectral discretisation in space followed by exponential
splittings for time-propagation \cite{bao02ots}. Methods based on Lanczos
iterations \cite{saad92s}, which are particularly effective for problems in the
atomic scaling ($\ve=1$), are also found to be ineffective in the
semiclassical regime, where the exponent involved is of huge spectral size,
scaling as $\O{\ve^{-1}}$, see~\cite{hoclub97,lubich08}.
\rev{Effective
methods for highly oscillatory problems, which have a Hamiltonian structure
and are periodic in time, have been proposed in \cite{brugnano18b}
and are referred to as \emph{Hamiltonian boundary value methods}.
They have been applied in the context of Schr{\"o}dinger equations
in \cite{brugnano18c}.
}

In this regime the symmetric Zassenhaus splittings of \cite{bader14eaf} are
found to be very effective. These are asymptotic exponential splittings where
 exponents scale with powers of the small parameter $\ve$ and,
consequently, become progressively small. The small size of these exponents
allows a very effective approximation via Lanczos iterations despite
reasonably large time steps~\cite{bader14eaf}.


Needless to say, the oscillations in the solution change in time. Higher
oscillations require smaller time steps, which increases the computational cost.
Thus, we face the usual tradeoff between two
competing concerns -- smaller time steps for higher accuracy and larger time
steps for lower cost. Time adaptivity has been a well developed approach to
arrive at the optimal compromise by keeping the time steps as large as
possible, so long as a prescribed error tolerance is not exceeded.
State-of-the-art step-size choice is based on firm theoretical ground by recent investigations of
digital filters from signal processing and control theory
\cite{soederlind01,soederlind03,soederlind06a}, which have also been demonstrated to
enhance computational stability \cite{soederlind06b}.
The advantages of adaptive selection of the time steps in the context of
splitting methods have been demonstrated for nonlinear Schr{\"o}dinger equations
in \cite{koarl2paper}, and for parabolic equations in \cite{quelletal15}.
It is found that in addition to a potential increase in the computational efficiency,
the reliability of the computations is enhanced. In an adaptive procedure,
the step-sizes are not guessed a priori but chosen automatically as
mandated by the smoothness of the solution, and the solution accuracy
can be guaranteed.

The aim of this paper is thus to design a time-adaptive method for the \schr
equation in the semiclassical regime. This is done by utilising the defect
based time adaptivity approach of \cite{auzingeretal13a,auzingeretal18a} for
the symmetric Zassenhaus splittings of \cite{bader14eaf}. These time
adaptivity schemes, which are very effective for classical splittings, turn
out to be successful for the asymptotic splittings of \cite{bader14eaf}
as well.

In Section~\ref{sec:def-ad} we describe the defect-based time-adaptive
approach. In Section~\ref{sec:symm-zass-spl} we present a variant of the
high-order Zassenhaus splittings of \cite{bader14eaf}, which turns out to be
more conducive in a time-adaptive approach. The practical algorithm used for
the estimation of the local error for this highly efficient method is derived
in Section~\ref{sub:localerror}. In Section~\ref{sec:num} we present
numerical experiments which demonstrate the efficacy of the time-adaptive
scheme.

\section{Defect-based time adaptivity} \label{sec:def-ad}
Time-adaptivity in numerical solutions of ODEs and PDEs
involves adjusting time step of a numerical scheme in order to keep the local
error in a single step below a specified error tolerance, ${\tt tol}$. The
local error in a single step of a one-step integrator with step-size $h$,
starting (without loss of generality) from $ \psi_0 $ at $t=0$,
\[ \psi_1 = \SSS(h) \psi_0, \]
where $\SSS(h)$ is the numerical evolution operator, is given by
\[ \LLL(h)\psi_0 = (\SSS(h)-\EEE(h))\psi_0, \]
where
\[ \psi(h) = \EEE(h) \psi_0\]
is the exact solution. Since the exact evolution operator for (\ref{eq:TDSE}),
\begin{equation*}
\EEE(h) = \ee^{\ii\,\dt (\ve \dx^2 - \ve^{-1} V)},
\end{equation*}
and the exact solution are not available, practical time-adaptivity methods
rely on accurate a~posteriori estimates of the local error, $\tilde{\LLL}(h)
\psi_0$, that can be computed along with the numerical solution. In this
manuscript, we will focus on defect-based estimators of the local error,
which are of the form
\begin{equation} \label{eq:defect-based-lerr-est}
{\widetilde\LLL}(h)\psi_0 = \frac{h}{p+1} \DDD(h)\psi_0 = \O{h^{p+1}}.
\end{equation}
Here, $ \DDD(h)\psi_0 $ is a computable defect term, measuring the local
quality of the approximation delivered by the numerical solution, $
\SSS(h)\psi_0 $. We will consider two different versions of the defect $
\DDD(h)\psi_0 $, see Sections~\ref{sec:dc} and~\ref{sec:ds}.

Once an estimate of the local error via (\ref{eq:defect-based-lerr-est}) is
available, a new step-size can then be chosen as
\[ h_{\mathrm{new}} =
(1-\alpha)\,h_{\mathrm{old}}
\sqrt[\leftroot{-2}\uproot{2}p+1]{\frac{{\tt tol}}{{\widetilde\LLL}(h)\psi_0}},
\quad \alpha > 0,\] where ${\tt tol}$ is the local error tolerance and
$p+1$ is the local order of the numerical scheme, $\SSS(h)=\EEE(h) +
\O{h^{p+1}}$. The factor of $(1-\alpha)$ ensures that we are more
conservative with the time step, always taking a slightly smaller time step
than predicted ($\alpha$ is usually taken to be a small positive number such
as $0.1$). When the local error estimate in a step
$\widetilde{\LLL}(h)\psi_0$ exceeds the error tolerance (${\tt tol}$), the
numerical propagation is run again with the smaller step
$h=h_{\mathrm{new}}$, otherwise the new time step is used for the next step.

\subsection{Classical defect-based estimator} \label{sec:dc}
We briefly recall the idea underlying~(\ref{eq:defect-based-lerr-est});
see for instance~\cite{auzingeretal13a}.
Thinking of the step-size $ h $ as a continuous variable and denoting it by $ t $
(somewhat more intuitive), the discrete flow $ \SSS(t)\psi_0 $ is a well-defined, smooth function of $ t $.
We call
\begin{equation}
\label{eq:defect-classical}
\DDD_c(t) := \partial_t \SSS(t) - H \SSS(t)
\end{equation}
the {\em classical defect}\/ (operator) associated with $ \SSS $, obtained by
plugging in the numerical flow into the given evolution
equation~(\ref{eq:TDSE}) (which is satisfied exactly by its exact flow $
\EEE(t) $).
Due to the definition of $ \DDD_c $, the local error enjoys the integral representation
\begin{subequations}
\begin{equation}\label{eq:locerr-int-c}
\LLL(h)\psi_0 = \int_0^h \EEE(h-t) \DDD_c(t)\psi_0\,dt = \O{h^{p+1}}.
\end{equation}
As argued in~\cite{auzingeretal13a}, this can be approximated in an asymptotically
correct way via its classical defect $ \DDD_c(h) \psi_0 $:
\begin{equation}\label{eq:locerr-appr-c}
{\widetilde\LLL}_c(h)\psi_0 := \frac{h}{p+1}\,\DDD_c(h)\psi_0 = \LLL(h)\psi_0 + \O{h^{p+2}} \quad \text{for}~ h \to 0.
\end{equation}
\end{subequations}
Evaluation of $ {\widetilde\LLL}_c(h)\psi_0 $
requires a single evaluation of the defect at $ t=h $,
the step-size actually used in the computation of $ \psi_1 = \SSS(h)\psi_0 $.\footnote{How to compute
 the defect at $ t=h $, for the given step-size $ h $ used in the computation,
 in the context of Zassenhaus splitting will be explained
 in~Section~\ref{sec:locerr-zass-split}.}

\subsection{Symmetrized defect-based estimator}  \label{sec:ds}
An alternative way of defining the defect was introduced in~\cite{auzingeretal18a,auzingeretal18c}.
It is based on the fact that the exact evolution operator
$ \EEE(t) $ commutes with the Hamiltionian $ H $, whence
\[
\partial_t \EEE(t)\psi_0 = H \EEE(t)\psi_0 = \EEE(t) H \psi_0
                         = \tfrac{1}{2}\{H,\EEE(t)\} \psi_0,
\]
with the anti-commutator $ \{ H,X \} = H X + X H $. This motivates the
definition of the {\em symmetrized defect}\/ (operator)
\begin{equation}\label{eq:defect-symmetrized}
\DDD_s(t):=
\partial_t \SSS(t) - \tfrac{1}{2}\{H,\SSS(t)\}.
\end{equation}
\begin{subequations}
Then, the local error $ \LLL $ enjoys the alternative integral representation
\begin{equation}\label{eq:locerr-int-s}
\LLL(h)\psi_0 = \int_0^h \EEE(\tfrac{h-t}{2}) \DDD_s(t) \EEE(\tfrac{h-t}{2})\psi_0\,\dd t,
\end{equation}
and by the same reasoning as for the classical defect,
this can be approximated in an asymptotically correct way via
the symmetrized defect $ \DDD_s(h)\psi_0 $\,:
\begin{equation}\label{eq:locerr-appr-s}
{\tilde\LLL}_s(h)\psi_0 := \frac{h}{p+1}\,\DDD_s(h)\psi_0
= \LLL(h)\psi_0 + \O{h^{p+2}}.
\end{equation}
\end{subequations}
Now suppose that $ \SSS $ is symmetric (time-reversible), i.e., it satisfies
$ \rev{\SSS(-t)\SSS(t)} = I $. Then, its order $ p $ is necessarily even, and moreover we even have
\begin{equation}\label{eq:locerr-appr-s-1}
{\widetilde\LLL}_s(h)\psi_0 = \LLL(h)\psi_0 + \O{h^{p+3}} \quad \text{for}~ h \to 0,
\end{equation}
see~\cite{auzingeretal18a,auzingeretal18c}. This means that in the symmetric case the deviation
$ ({\widetilde\LLL}_s(h) - \LLL(h))\psi_0 $ of the
local error estimate based on the symmetrized defect is of a better quality,
asymptotically for $ h \to 0 $, than
that one based on the classical defect.\footnote{Note that $ \DDD_c(h) = \O{h^{p}} $,
                                     $ \DDD_s(h) = \O{h^{p}} $, and
                                     $ \DDD_s(h) - \DDD_c(h) = \O{h^{p+1}} $.}
Moreover, evaluation of $ \DDD_s(h)\psi_0 $ is
typically only slightly more expensive compared to $ \DDD_c(h)\psi_0 $,
see~\cite{auzingeretal18a,auzingeretal18c} for the case
of conventional splittings and, in particular,
Section~\ref{sec:locerr-zass-split} below.

\section{Symmetric Zassenhaus splittings} \label{sec:symm-zass-spl}
Symmetric Zassenhaus splittings are asymptotic splittings introduced
in~\cite{bader14eaf} for solving the \schr equation in the semiclassical
regime~(\ref{eq:TDSE}). These splittings are derived by working in infinite
dimensional space, prior to spatial discretisation, which enables a more
effective exploitation of commutators, since some of them vanish, while the
rest end up being much smaller (in a sense of spectral radius) than
generically expected.

To manage the $\O{\ve}$ wavelength \rev{oscillations} in space and time, which arise
due to the presence of the small semiclassical parameter $\ve$, the following
relations have been established to be useful for the spatio-temporal
resolution in these schemes.
\begin{equation}
\label{eq:eps-h-relations}
\Delta x = \O{\ve}, \qquad \dt = \O{\ve^{\sigma}},\ \sigma \in (0,1].
\end{equation}

The Zassenhaus splitting that we will consider in this paper is
\begin{subequations}
\begin{equation}
\label{eq:Zsigma}
\SSS(h) =
\ee^{\frac12 W^{[0]}}
 \ee^{\frac12 W^{[1]}}
  \ee^{\frac12 W^{[2]}}
   \ee^{\WWW^{[3]}}
    \ee^{\frac12 W^{[2]}}
     \ee^{\frac12 W^{[1]}}
      \ee^{\frac12 W^{[0]}}
= \EEE(h) + \O{\ve^{7\sigma - 1}},
\end{equation}
where
\begin{eqnarray}
  \label{eq:Wzass1-modified}
  W^{[0]}&=&\ii\,\dt\,\ve\,\dx^2=\O{\ve^{\sigma-1}}, \\
  \notag
  W^{[1]}&=&-\ii\,\dt\,\ve^{-1} V=\O{\ve^{\sigma-1}}, \\
  \notag
  W^{[2]}&=&\tfrac16\ii\,\dt^3\ve^{-1}(\dx V)^2-\tfrac{1}{24}\ii\,\dt^3\ve(\dx^4V)+\tfrac{1}{6}\ii\,\dt^3\ve \ang{\dx^2V}{2}=\O{\ve^{3\sigma-1}}, \\
  \notag
  \WWW^{[3]}&=& -\tfrac{7}{120}\ii\,\dt^5\ve^{-1} (\dx^2V)(\dx V)^2 +\tfrac{1}{30}\ii\,\dt^5\ve \ang{(\dx^2V)^2  - 2(\dx^3V)(\dx V)}{2} \notag \\
  &&~~{} - \tfrac{1}{120}\ii\,\dt^5\ve^3
  \ang{\dx^4V}{4}=\O{\ve^{5\sigma-1}}.
   \notag
\end{eqnarray}
\end{subequations}

\noindent Here,
\begin{equation} \label{eq:<f>k}
\ang{f}{k} := \tfrac12 \left( f \circ \dx^k + \dx^k \circ f \right),
\quad k \geq 0 \qquad (f \in \mathrm{C}_p^{\infty}(I;\BB{R}))
\end{equation}
are the symmetrized differential operators which first appeared
in~\cite{bader14eaf} and have been studied in detail in
\cite{psalgebra,Singh17}. These operators are used for preservation of
stability under discretization after simplification of commutators.

Once a differential operator such as $\dx^2$ is discretised to a symmetric
differentiation matrix $\K_2$ via spectral collocation, its spectral size
grows as $\Delta x = \O{\ve}$ decreases, since
\[ \norm{2}{\K_k} = \O{(\Delta x)^{-k}} = \O{\ve^{-k}}.\]
Keeping eventual discretisation with the scaling (\ref{eq:eps-h-relations})
in mind, we use a shorthand $\dx^k = \O{\ve^{-k}}$ for the undiscretised and
unbounded operator as well. The symmetrized differential operator
$\ang{f}{k}$ discretises to the form
\[ \widetilde{\ang{f}{k}} = \tfrac12 \left( \DDD_f \K^k + \K^k \DDD_f \right), \]
where $\DDD_f$ is a diagonal matrix with values of $f$ at the grid points.
Its spectral radius grows as $\O{\ve^{-k}}$, assuming that $f$ is independent
of $\ve$, and we abuse notations once more to say $\ang{f}{k} =
\O{\ve^{-k}}$, for short.

The splitting (\ref{eq:Zsigma}) possesses several favorable features. First of
all, the critical quantities like time step $h$, spatial step $\Delta x$ and
semiclassical parameter $\ve$ are tied together by
relation~(\ref{eq:eps-h-relations}). As a consequence, the splitting error is
expressed via a universal quantity $\mathcal{O}(\ve^{7\sigma-1})$, and the
error constant does not hide any critical quantities.

The exponentials involved in the splitting are easily computable after
spatial discretization. For example, spectral collocation transforms
$W^{[0]}$ into a symmetric circulant matrix, whence its exponential may be
computed by the Fast Fourier Transform. The diagonal matrix  $W^{[1]}$ can be
exponentiated directly. Neither $W^{[2]}$ nor $\mathcal{W}^{[3]}$ are
structured  favourably, however their spectral radius is small enough to
exponentiate them with a small number (say, 3 or 4) of Lanczos iterations.

\begin{remark}
Since $h=\O{\ve^\sigma}$ the splitting {\rm (\ref{eq:Zsigma})} which features an error
of $\O{\ve^{7\sigma-1}}$ could loosely be considered a sixth-order splitting
with an error constant scaling as $\ve^{-1}$. This interpretation is not
strictly correct, however, since the error estimate holds for {$\sigma \in
(0,1]$} and $\ve \rightarrow 0$. For a finite $\ve = \ve_0$ and $h
\rightarrow 0$ (i.e. $\sigma \rightarrow \infty$), the asymptotic Zassenhaus
splitting {\rm (\ref{eq:Zsigma})} is a highly efficient fourth-order method with a very
small error constant. This disparity in the different asymptotic limits
occurs because the derivation of {\rm (\ref{eq:Zsigma})} involves discarding terms of
size $\O{h^5 \ve}$, for instance. These terms are $\O{\ve^{7 \sigma -1}}$
under the asymptotic scaling {\rm (\ref{eq:eps-h-relations})}, $\ve \rightarrow 0, h =
\O{\ve^\sigma}, \sigma \leq 1$, but are $\O{h^5}$ under a fixed $\ve=\ve_0$
and $h \rightarrow 0$. The consequence for the time-adaptivity algorithm
{\rm (\ref{eq:defect-based-lerr-est})} is that we must use $p=4$, not $p=6$.
\end{remark}

\begin{remark}
Many different variants of this splitting are possible (including, but not
limited to the case when we start with $W^{[0]} = -\ii\,\dt\,\ve^{-1} V$).
Application of these variants can just as easily be considered, however will
not be the focus of the present paper.
\end{remark}

\begin{remark}
The exponents {\rm (\ref{eq:Wzass1-modified})} differ from the splitting
described in {\rm \cite{bader14eaf}} in one minor aspect -- \rev{the
 term $-\tfrac{1}{24}\ii\,\dt^3\ve(\dx^4V)$} has been moved from
$\WWW^{[3]}$ to $W^{[2]}$. \rev{This term is of size $\O{\ve^{3\sigma+1}}$,
which is smaller than $\O{\ve^{3\sigma-1}}$ and thus may also be combined
with $W^{[2]} = \O{\ve^{3\sigma-1}}$. In {\rm \cite{bader14eaf}} this term is
combined with $\WWW^{[3]} = \O{\ve^{5\sigma-1}}$ working under the assumption
$\sigma \leq 1$ (time steps larger than $\O{\ve}$) since
$\O{\ve^{3\sigma+1}}$ is also smaller than $\O{\ve^{5\sigma-1}}$ under
$\sigma \leq 1$.} Consequently, all exponents feature a single power of
$\dt$. This change makes no difference to the overall order of the scheme but
makes the computation of defect easier.
\end{remark}

\section{Local error estimator for Zassenhaus splitting} \label{sec:locerr-zass-split}
\label{sub:localerror}
\begin{remark}
In the sequel, the argument $ (t) $ or $ (h) $, respectively, is suppressed
whenever the meaning is obvious.
\end{remark}

\subsection{Time derivative of the discrete evolution operator}
\label{sub:S-time-derivative} The time-adaptivity algorithm is based on the
estimation of the local error in each step, which in turn requires the
computation of the defect at discrete time $t=h$, see Section~\ref{sec:def-ad}.
This requires the computation of $ \partial_t \SSS(h) = \partial_t
\SSS(t)|_{t=h} $ For splitting schemes such as~(\ref{eq:Zsigma}), this boils down
to an application of the product rule, involving computation of the time
derivative of each individual exponential.

For exponentials such as
\begin{subequations}
\begin{equation}
\label{eq:S}
\SSS_j(t) = \exp(t^{n_j} \RRR^{[j]}),
\end{equation}
where $\RRR^{[j]}$ is independent of $t$, the time derivative is
\begin{equation}
\label{eq:dtS}
\partial_t \SSS_j(t) = n_j t^{n_j-1} \RRR^{[j]} \exp(t^{n_j} \RRR^{[j]}).
\end{equation}
\end{subequations}
Now, consider the (modified) Zassenhaus
splitting~(\ref{eq:Zsigma})\,\&\,(\ref{eq:Wzass1-modified}) for~(\ref{eq:TDSE}), with the
discrete evolution operator re-written in the form
\begin{eqnarray*}
\notag \SSS(t) &=& \SSS_0\,\SSS_1\,\SSS_2\,\SSS_3\,\SSS_2\,\SSS_1\,\SSS_0 \\
\label{eq:Zt}&=& \ee^{t^{n_0} \RRR^{[0]}} \ee^{t^{n_1} \RRR^{[1]}}
\ee^{t^{n_2} \RRR^{[2]}} \ee^{t^{n_3} \RRR^{[3]}} \ee^{t^{n_2} \RRR^{[2]}}
\ee^{t^{n_1} \RRR^{[1]}} \ee^{t^{n_0} \RRR^{[0]}},
\end{eqnarray*}
where $\SSS_j = \SSS_j(t) $, $ n_0 = 1,\ n_1= 1,\ n_2 = 3,\ n_3 = 5 $, and
\begin{eqnarray}
  \label{eq:Rzass1}
  \RRR^{[0]}&=&\tfrac12 \ii\,\ve\,\dx^2, \\
  \notag
  \RRR^{[1]}&=&-\tfrac12 \ii\,\ve^{-1} V, \\
  \notag
  \RRR^{[2]}&=&\tfrac{1}{12}\ii\,\ve^{-1}(\dx V)^2-\tfrac{1}{48}\ii\,\ve(\dx^4V)+\tfrac{1}{12}\ii\,\ve \ang{\dx^2V}{2}, \\
  \notag
  \RRR^{[3]}&=& -\tfrac{7}{120}\ii\,\ve^{-1} (\dx^2V)(\dx V)^2 +\tfrac{1}{30}\ii\,\ve \ang{(\dx^2V)^2  - 2(\dx^3V)(\dx V)}{2}  -\tfrac{1}{120}\ii\,\ve^3\ang{\dx^4V}{4}.
  \notag
\end{eqnarray}
The derivative of the flow for a splitting scheme can be expressed via the
product rule,
\begin{eqnarray*}
\partial_t \SSS(t)
&=& n_0t^{n_0-1} \SSS_0\,\SSS_1\,\SSS_2\,\SSS_3\,\SSS_2\,\SSS_1\,\RRR^{[0]}\,\SSS_0\, \\
&& {} + n_1t^{n_1-1} \SSS_0\,\SSS_1\,\SSS_2\,\SSS_3\,\SSS_2\,\SSS_1\,\RRR^{[1]}\,\SSS_0\, \\
&& {} + n_2t^{n_2-1} \SSS_0\,\SSS_1\,\SSS_2\,\SSS_3\,\RRR^{[2]}\,\SSS_2\,\SSS_1\,\SSS_0\, \\
&& {} + n_3t^{n_3-1} \SSS_0\,\SSS_1\,\SSS_2\,\SSS_3\,\RRR^{[3]}\,\SSS_2\,\SSS_1\,\SSS_0\, \\
&& {} + n_2t^{n_2-1} \SSS_0\,\SSS_1\,\RRR^{[2]}\,\SSS_2\,\SSS_3\,\SSS_2\,\SSS_1\,\SSS_0\, \\
&& {} + n_1t^{n_1-1} \SSS_0\,\SSS_1\,\RRR^{[1]}\,\SSS_2\,\SSS_3\,\SSS_2\,\SSS_1\,\SSS_0\, \\
&& {} + n_0t^{n_0-1} \RRR^{[0]}\,\SSS_0\,\SSS_1\,\SSS_2\,\SSS_3\,\SSS_2\,\SSS_1\,\SSS_0\, \\
&=&  \SSS_0\,\SSS_1\,\SSS_2\,\SSS_3\,\SSS_2\,\SSS_1\,(n_0t^{n_0-1} \RRR^{[0]} + n_1t^{n_1-1} \RRR^{[1]})\,\SSS_0\, \\
&& {} + \SSS_0\,\SSS_1\,\SSS_2\,\SSS_3\,(n_2t^{n_2-1} \RRR^{[2]} + n_3t^{n_3-1} \RRR^{[3]})\,\SSS_2\, \SSS_1\, \SSS_0\, \\
&& {} + \SSS_0\,\SSS_1\,(n_2t^{n_2-1} \RRR^{[2]} +  n_1t^{n_1-1} \RRR^{[1]})\,\SSS_2\,\SSS_3\, \SSS_2\, \SSS_1\, \SSS_0\, \\
&& {} + n_0t^{n_0-1} \RRR^{[0]} \SSS_0\,\SSS_1\,\SSS_2\,\SSS_3\,\SSS_2\,\SSS_1\,\SSS_0\, \\
&=&  \underline{\SSS_0\,\SSS_1} \left\{\underline{\SSS_2\,\SSS_3}
\left[\underline{\SSS_2\,\SSS_1}(n_0t^{n_0-1} \RRR^{[0]} + n_1
t^{n_1-1}\RRR^{[1]})\,\SSS_0\,\right.\right.\\
&&\left. +(n_2t^{n_2-1} \RRR^{[2]} + n_3t^{n_3-1} \RRR^{[3]})\,\SSS_2\, \SSS_1\, \SSS_0\,\right]\\
&& {} +\left. (n_2t^{n_2-1} \RRR^{[2]} +  n_1t^{n_1-1} \RRR^{[1]})\,\SSS_2\,\SSS_3\, \SSS_2\, \SSS_1\, \SSS_0\,\right\}\\
&& {} + n_0t^{n_0-1} \RRR^{[0]}
\SSS_0\,\SSS_1\,\SSS_2\,\SSS_3\,\SSS_2\,\SSS_1\,\SSS_0.
\end{eqnarray*}
Here the \underline{underlined} exponentials are the ones that need to be
computed freshly. The rest can be dealt with by storing intermediate values
from evaluation of the splitting scheme: In order to compute $ \partial_t
\SSS(h)$, we store
\begin{eqnarray}
\label{eq:v}
v_1 &=& \SSS_0\,\psi_0, \\
\notag v_2 &=& \SSS_2\,\SSS_1\,\SSS_0\,\psi_0 = \SSS_2\,\SSS_1 v_1, \\
\notag v_3 &=& \SSS_2\,\SSS_3\,\SSS_2\,\SSS_1\,\SSS_0\,\psi_0 = \SSS_2\,\SSS_3 v_2, \\
\notag (\psi_1 =\,)~v_4 &=& \SSS_0\,\SSS_1\,\SSS_2\,\SSS_3\,\SSS_2\,\SSS_1\,\SSS_0\,\psi_0 = \SSS_0\,\SSS_1 v_4,
\end{eqnarray}
which are anyway computed during evaluation of the numerical scheme.
In addition, we compute and store
\begin{eqnarray}
\label{eq:w}
w_1 &=& (n_0 h^{n_0-1} \RRR^{[0]} + n_1 h^{n_1-1}\RRR^{[1]}) v_1 = (\RRR^{[0]} + \RRR^{[1]}) v_1, \\
\notag w_2 &=& (n_2 h^{n_2-1} \RRR^{[2]} + n_3 h^{n_3-1} \RRR^{[3]}) v_2 = (3 h^2 \RRR^{[2]} + 5 h^4 \RRR^{[3]}) v_2, \\
\notag w_3 &=& (n_2 h^{n_2-1} \RRR^{[2]} + n_1 h^{n_1-1} \RRR^{[1]}) v_3 = (3 h^2 \RRR^{[2]} + \RRR^{[1]}) v_3, \\
\notag w_4 &=& n_0 h^{n_0-1} \RRR^{[0]} v_4 = \RRR^{[0]} v_4.
\end{eqnarray}
Then,
\begin{equation}
\label{eq:partialS(h)psi0}
\partial_t \SSS(h)\psi_0 =
\SSS_0(h)\,\SSS_1(h)
 \left\{\SSS_2(h)\,\SSS_3(h)
  \left[\SSS_2(h)\,\SSS_1(h) w_1 + w_2
  \right] + w_3
 \right\} + w_4.
\end{equation}
Thus, we need to compute six \rev{exponentials appearing in
(\ref{eq:partialS(h)psi0}) in addition to the seven exponentials required in
the Zassenhaus splitting (\ref{eq:Zsigma})}. For a scheme of order six, we
need a total of 13 exponentials.

\rev{
\begin{remark}
For the order four method
\[ \ee^{\frac12 W^{[0]}}
 \ee^{\frac12 W^{[1]}}
  \ee^{W^{[2]}}
     \ee^{\frac12 W^{[1]}}
      \ee^{\frac12 W^{[0]}}
= \EEE(h) + \O{\ve^{5\sigma - 1}},\] five exponentials are required for the
Zassenhaus splitting and we have verified that four additional exponentials
are required for time adaptivity, making for a total of nine exponentials.
\end{remark}}

\begin{remark}
In a practical, memory-efficient implementation, the $ w_j $ are computed `on
the fly' along with the $ v_j $, via alternating updates of arrays $ v $ and
$ w $.
\end{remark}

\begin{remark}
In the above computations, the only feature specific to the Zassenhaus
splitting is the use of~(\ref{eq:dtS}) for the derivative of the
exponential. Thus the derivation and the combination of exponentials
works in an analogous way for other splittings, with different time derivatives.
\end{remark}

\subsection{Practical evaluation of the classical defect~{\rm (\ref{eq:defect-classical})}}
We need to compute the classical defect at discrete time $ t=h $,
\begin{equation}
\label{eq:defectdiscrete-classical}
\DDD_c(h)\psi_0 = \partial_t \SSS(h)\psi_0 - H \psi_1.
\end{equation}
Here, $ \psi_1 = \SSS(h) \psi_0 $,
$H = 2 (\RRR^{[0]}+\RRR^{[1]})$ is the Hamiltonian,
and $ \partial_t \SSS(h)\psi_0 $ is the time derivative
according to Section~\ref{sub:S-time-derivative}, see~(\ref{eq:partialS(h)psi0}).

In addition, a small amount of work can be saved by exploiting
the relation $ w_4 = \RRR^{[0]} v_4 $ (see~\ref{eq:w}). Then,
computation of $ \DDD_c(h)\psi_0 $ amounts to the evaluation of
\begin{equation*}
\DDD_c(h)\psi_0 = \SSS_0(h)\,\SSS_1(h)
                \left\{\SSS_2(h)\,\SSS_3(h)
                \left[\SSS_2(h)\,\SSS_1(h)w_1 + w_2\right] + w_3\right\}
- (\RRR^{[0]} + 2\RRR^{[1]}) v_4.
\end{equation*}

\subsection{Practical evaluation of the symmetrized defect~{\rm (\ref{eq:defect-symmetrized})}}
Instead of~(\ref{eq:defectdiscrete-classical}), we need to compute
\begin{equation}
\label{eq:defectdiscrete-symmetrized}
\DDD_s(h)\psi_0
= \partial_t \SSS(h)\psi_0
- \tfrac12 \big( H \psi_1 + \SSS(h) H \psi_0 \big).
\end{equation}
With~(\ref{eq:partialS(h)psi0}), we have
(again making use of $ w_4 = \RRR^{[0]} v_4 $)
\begin{eqnarray*}
\DDD_s(t)\,\psi_0 &=& \nomark{\SSS_0\,\SSS_1\,}
\left\{\nomark{\SSS_2\,\SSS_3\,} \left[\nomark{\SSS_2\,\SSS_1\,} w_1 + w_2\right] + w_3\right\} +w_4
- \tfrac12 \SSS H \psi_0 - \tfrac12 H \SSS \psi_0 \\
&=& \nomark{\SSS_0\,\SSS_1\,}
\left\{\nomark{\SSS_2\,\SSS_3\,} \left[\nomark{\SSS_2\,\SSS_1\,} w_1 + w_2\right] + w_3\right\} +w_4 \\
&& {} - \tfrac12 \SSS_0\,\SSS_1\,\SSS_2\,\SSS_3\,\SSS_2\,\SSS_1\,\SSS_0\,H \psi_0 - \tfrac12 H \SSS \psi_0 \\
&=& \nomark{\SSS_0\,\SSS_1\,} \left\{\nomark{\SSS_2\,\SSS_3\,}
\left[\nomark{\SSS_2\,\SSS_1\,} (w_1 - \tfrac12
\underline{\SSS_0\,} H \psi_0) + w_2\right] + w_3\right\} + w_4 - \tfrac12 H \SSS \psi_0 \\
&=&\nomark{\SSS_0\,\SSS_1\,} \left\{\nomark{\SSS_2\,\SSS_3\,}
\left[\nomark{\SSS_2\,\SSS_1\,} (w_1 - \tfrac12
\underline{\SSS_0\,} H \psi_0) + w_2\right] + w_3\right\} + w_4 - (\RRR^{[0]}+ \RRR^{[1]}) v_4\\
 &=& \nomark{\SSS_0\,\SSS_1\,} \left\{\nomark{\SSS_2\,\SSS_3\,}
\left[\nomark{\SSS_2\,\SSS_1\,} (w_1 - \tfrac12 \underline{\SSS_0\,} H \psi_0) +
w_2\right] + w_3\right\} - \RRR^{[1]} v_4\\
 &=& \nomark{\SSS_0\,\SSS_1\,} \left\{\nomark{\SSS_2\,\SSS_3\,}
\left[\nomark{\SSS_2\,\SSS_1\,} (w_1 - \underline{\SSS_0\,} (\RRR^{[0]} +
\RRR^{[1]}) \psi_0) + w_2\right] + w_3\right\} - \RRR^{[1]} v_4.
\end{eqnarray*}
Here the \underline{underlined} exponential is the one that needs to be
computed in addition, compared to evaluation  of $ \DDD_c(h)\psi_0 $. Of
course, this additional effort is marginal.

Thus, computation of the symmetrized defect at discrete time $ t=h $ amounts to
evaluation of
\begin{align*}
& \DDD_s(h)\psi_0 = \\
& \SSS_0(h)\,\SSS_1(h)\{\SSS_2(h)\,\SSS_3(h)\,
[\SSS_2(h)\,\SSS_1(h) (w_1 - \SSS_0(h)(\RRR^{[0]} +
\RRR^{[1]}) \psi_0) + w_2] + w_3\} - \RRR^{[1]} v_4.
\end{align*}

\section{Numerical experiments} \label{sec:num}
For our numerical experiments, we will use the wave-packets
\begin{equation*}
\label{eq:Lu0}
\psi_L(x) = \varphi\left(x;\tfrac{\ve}{4},-\tfrac{3}{4},\tfrac{1}{10}\right), \quad \psi_M(x) = \varphi\left(x;\ve,\tfrac{9}{2},0\right),
\end{equation*}
as initial conditions for two experiments. Here
\begin{equation*}
\label{eq:wavepacket}
\varphi(x;\delta,x_0,k_0) =
(\delta \pi)^{-1/4} \exp\left(\ii k_0 \frac{(x-x_0)}{\delta}-\frac{(x-x_0)^2}{2 \delta}\right),
\end{equation*}
is a wave-packet with a spread of $\delta$, mean position $x_0$ and mean momentum $k_0$.
We will consider the evolution of $\psi_L$ as it heads towards the lattice potential $V_L$ and the evolution of $\psi_M$ as it oscillates in the Morse potential $V_M$, which are given by
\begin{equation*}
\label{eq:VL}
V_L(x) = \rho(4 x - 1)\sin(20\pi x) + \tfrac{1}{10}\rho(x/5)\sin(4\pi x), \quad V_M(x) = (1 - \ee^{-(x - 5)/2})^2,
\end{equation*}
respectively,
where
\[
\rho(x) =
\begin{cases}
\,\exp\left(-1/(1-x^2)\right) &\text{for } |x| < 1, \\
\,0 &\text{otherwise, }
\end{cases}\]
is a bump function.

\begin{figure}[ht]
       \centering
       \includegraphics[width=0.45\textwidth]{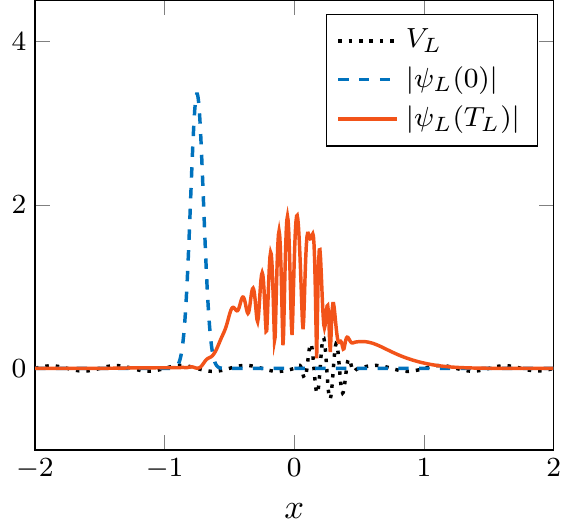}
       \includegraphics[width=0.459\textwidth]{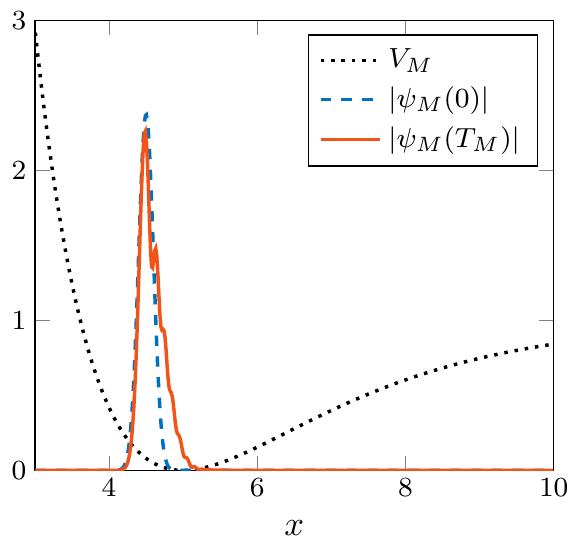}

    \caption{Initial wave-packet $\psi_L(0)$ evolves to the
    final wave-packet $\psi_L(T_L)$ at time $T_L=1$ under $V_L$ and $\ve=10^{-2}$ {\em (left)} and $\psi_M(0)$ evolves to $\psi_M(T_M)$ at $T_M=20$ under $V_M$ and $\ve=10^{-2}$ {\em (right)}.}
    \label{fig:experiments}
\end{figure}

We consider the behaviour at the final time $T_L=1$ for the first experiment
and $T_M=20$ for the second experiment under the choices of $\ve=10^{-2},
10^{-3}$ and $10^{-4}$. The spatial domain is chosen as $[-2,2]$ for the
lattice potential and $[3,10]$ for the Morse potential, and we impose
periodic boundary conditions. The evolution of these wavepackets under the
choice $\ve = 10^{-2}$ is shown in Fig.~\ref{fig:experiments}.

\begin{figure}[ht]
       \centering
       \includegraphics[width=0.999\textwidth]{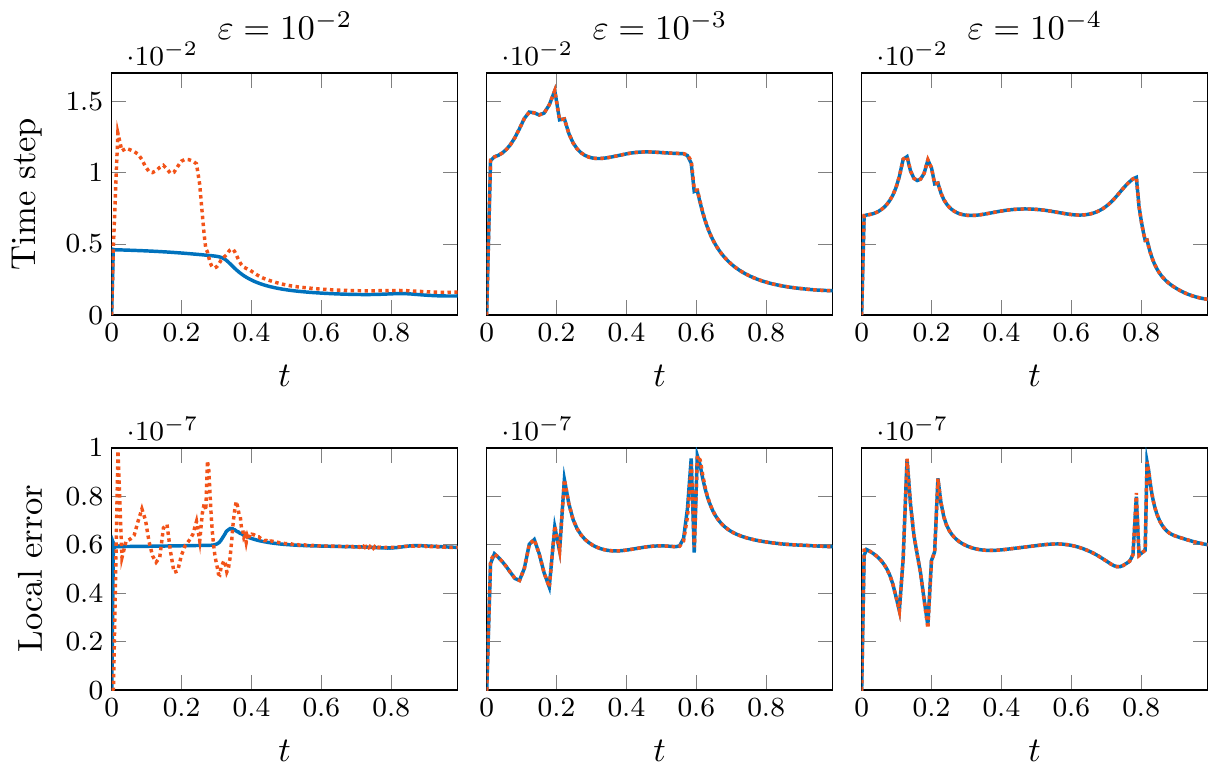}
    \caption{{\bf [Experiment 1: $\psi_L$ under $V_L$ for $t\in[0,1]$]} {\em (top row)} Time step  chosen by the adaptive procedure in consecutive steps of the numerical procedure,
    starting from a highly conservative initial guess of $h_0=10^{-9}$ and local error tolerance {\tt tol}$=10^{-7}$;
    {\em (bottom row)} Local error estimate in each step. The three columns show the behaviour in the regimes $\ve=10^{-2},10^{-3}$ and $10^{-4}$.
    The classical defect is depicted with solid blue line and the symmetric defect with dotted orange line.
    }
    \label{fig:adaptivityLattice}
\end{figure}

\begin{figure}[ht]
       \centering
       \includegraphics[width=0.999\textwidth]{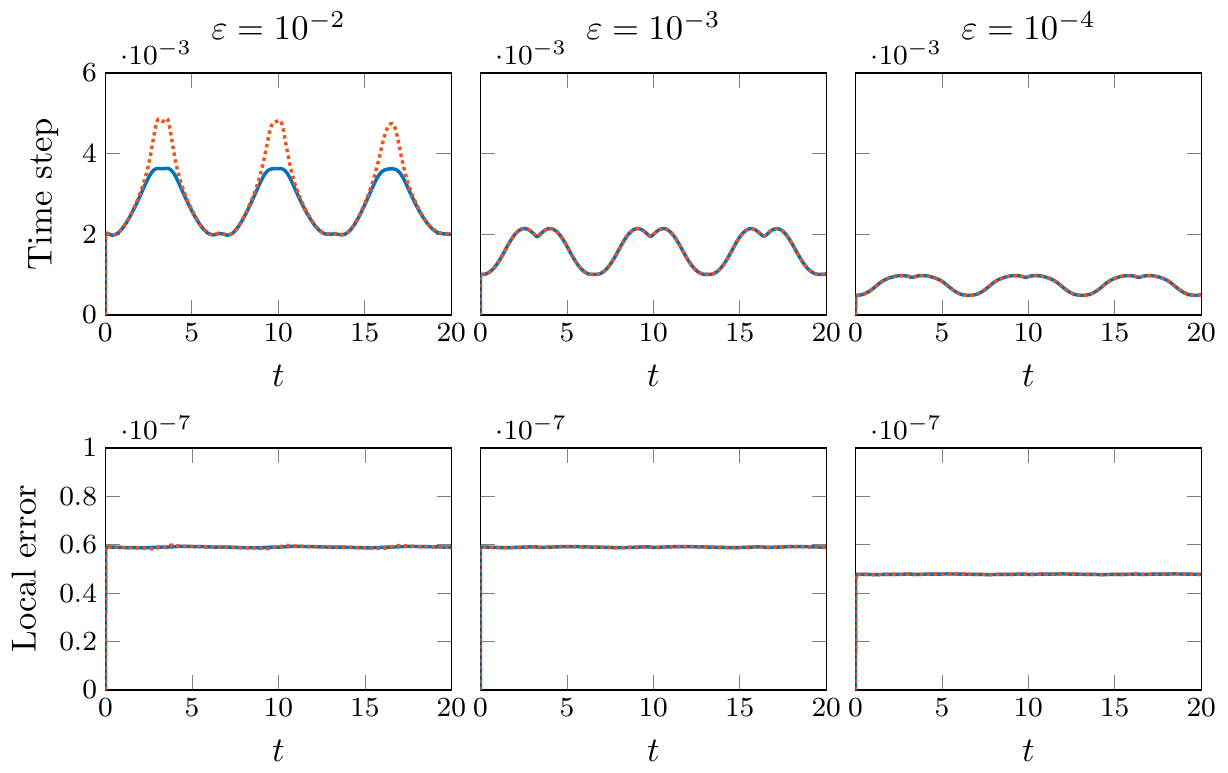}
    \caption{{\bf [Experiment 2: $\psi_M$ under $V_M$ at $T_M=20$]} {\em (top row)} Time step chosen by the adaptive procedure in consecutive steps of the numerical procedure,
    starting from a highly conservative initial guess of $h_0=10^{-9}$ and local error tolerance {\tt tol}$=10^{-7}$;
    {\em (bottom row)} Local error estimate in each step. The three columns show the behaviour in the regimes $\ve=10^{-2},10^{-3}$ and $10^{-4}$.
    The classical defect is depicted with solid blue line and the symmetric defect with dotted orange line.
    }
    \label{fig:adaptivityMorse}
\end{figure}

As one can observe in Figures~\ref{fig:adaptivityLattice} and \ref{fig:adaptivityMorse}, the procedure quickly
adapts the time step to ensure that the local accuracy is within the
specified threshold, which is taken to be ${\tt tol} = 10^{-7}$ here.

In the first
experiment, we use $M = 750, 1750$ and $15000$ equispaced points for spatial
discretization in the cases $\ve=10^{-2},10^{-3}$ and $10^{-4}$,
respectively, while $M=500,1500$ and $10000$ are used for the three choices of $\ve$ in the second experiment.

\begin{table}[!ht]
\centering
\begin{scriptsize}
\begin{tabular}{|c|c|c|c|c|c|c|c|}
  \cline{3-8}
  \multicolumn{2}{c|}{} & \multicolumn{3}{c|}{$\psi_L$ under lattice potential } & \multicolumn{3}{c|}{$\psi_M$ under Morse potential} \Tstrut\Bstrut\\
  \hline
       &   & Global  & Time & Run & Global  & Time & Run \Tstrut\\
  $\ve$ & h & error   & steps & time (s) & error   & steps & time (s) \Bstrut\\
       \hline
                    & \mbox{} adaptive (c)& $1.29 \times 10^{-6}$  & 490 & 19.5 & $2.69 \times 10^{-4}$  & 7772 & \ 320.0 \Tstrut\\
                    & adaptive (s) & $2.04 \times 10^{-5}$ & 394 & 18.7 & $2.94 \times 10^{-4}$ & 7448 & 316.0 \\
  $10^{-2}$   & smallest (c) & $1.96\times 10^{-7}$  & 746  & 11.9  & $1.92 \times 10^{-4}$  & 10102  & 191.5\\
                    & largest (c) & $5.53\times 10^{-4}$   & 217  & 4.1 & $6.43 \times 10^{-4}$   & 5502  & 104.1\\
                    & $h = \ve$ & $1.81\times 10^{-2}$  & 100 & 2.4 & $4.51 \times 10^{-3}$    & 2000 & 52.0 \Bstrut\\
  \hline
                      & adaptive (c) & $1.49 \times 10^{-7}$ & 216 & 12.9 & $9.19 \times 10^{-4}$   & 13249 & 955.1\Tstrut\\
                      & adaptive (s) & $1.53 \times 10^{-7}$ & 216 & 14.3 & $9.20 \times 10^{-4}$   & 13245 & 1021.5\\
  $10^{-3}$   & smallest (c)& $1.18\times 10^{-8}$   & 587 & 14.4 & $5.00 \times 10^{-4}$       & 19841 & 572.8\\
                    & largest (c)& $1.62\times 10^{-3}$  & 64 & 2.9 & $2.21 \times 10^{-3}$        & 9329 & 315.5  \\
                    & $h = \ve$ & $3.96\times 10^{-9}$   & 1000 & 20.4 & $4.92 \times 10^{-4}$         & 20000 & \ 601.4 \Bstrut\\
  \hline
                        & adaptive (c) & $9.70 \times 10^{-9}$  & 215 & 43.7 & $1.64 \times 10^{-3}$  & 27358 & 12029.7\Tstrut\\
                        & adaptive (s) & $9.75 \times 10^{-9}$  & 215 & 42.6 & $1.64 \times 10^{-3}$  & 27343 & 14237.4\\
  $10^{-4}$   & smallest (c) & $5.2\times 10^{-10}$   & 925 & 80.2 & $1.02\times 10^{-3}$   & 42555 & 11380.8\\
                    & largest (c) & $7.12\times 10^{-4}$  & 90 & 12.0 & $3.20 \times 10^{-3}$  & 21307 & 6435.6\\
                    & $h = \ve$ & $8.7\times 10^{-13}$  & 10000 & 659.8 & $ 4.95\times 10^{-5}$  & 200000 & 60755.0\Bstrut\\
  \hline
\end{tabular}
\end{scriptsize}
\caption{Results of different time-stepping strategies with an imposed tolerance of $10^{-7}$ for the lattice potential (left part of the table)
and the Morse potential (right), $\varepsilon\in\{10^{-2},10^{-3},10^{-4}\}.$ We compare the global errors with the
computational effort (number of time steps respectively computation time). Adaptive strategies based on
either the classical defect `(c)' or the symmetrized defect `(s)' are compared with fixed time steps chosen
as either the smallest adaptively determined step-size, the largest adaptive step-size or the asymptotic scaling $h=\varepsilon$.}
\label{tab:errors}
\end{table}

%
%
%

In Table~\ref{tab:errors}, we show the total number of time steps required
for maintaining the local accuracy of $10^{-7}$ for three cases: $\ve =
10^{-2}, 10^{-3}$ and $10^{-4}$. Also presented are the global $\mathrm{L}^2$
accuracies of these solutions at the final time ($T_L=1$ and $T_M=20$) and the computational time. These are
compared to the (non-adaptive) case when the time step is fixed at (i) the
finest value used by the adaptive method (ii) the coarsest value used by the
adaptive method (iii) $h=\ve$. The reference solutions for computing global
errors are produced by using {\sc Matlab}'s {\tt expm} function, while using
a much finer spatial grid.
The number of Lanczos iterations are chosen automatically to meet the
accuracy requirements using the a~priori error bounds of
\cite{hochbruck97ksa}.

We find that for small $\ve$, the asymptotic scaling $h=\ve$ yields the smallest error,
but the computational effort is prohibitive. The symmetrized error estimate invariantly
generates coarser time grids than the classical version, but the additional computational effort implies
that the picture is ambivalent in that the computational effort is smaller in only
about half of the cases. When the largest step generated by the adaptive strategy `(c)' is used
in a uniform grid, computations are fast but not sufficiently accurate, while use of the smallest
adaptive time step is prohibitively expensive in the case of the lattice potential
(for the Morse potential, the reduction in the number of time steps is outweighed
by the computational effort for the error estimate). Together, these two observations
imply that adaptive step-size choice constitutes the appropriate means to determine
grids that reproduce the solution in a reliable and efficient way, and excels over
a priori choice of the step-size. In particular, the guess of an optimal step-size
is commonly not feasible, while the adaptive strategy finds the time steps
automatically and guarantees the desired level of accuracy.

\section{Conclusions}\label{sec:concl}

We have investigated adaptive strategies used in conjunction with asymptotic
Zassenhaus splitting for the solution of linear time-dependent Schr{\"o}dinger
equations in the semiclassical regime. The local time-stepping is based
on defect-based error estimation strategies. We have demonstrated that
adaptivity provides a means to reliably and efficiently achieve a desired
level of accuracy especially for a small semiclassical parameter and thus
excels over fixed time steps in many cases. A major benefit is in addition that
no near-optimal eq\rev{u}idistant step-size has to be guessed a priori.




\end{document}